\documentclass[12pt]{article}
\usepackage{amsfonts, amsmath, amssymb}

\newcommand{\qed}{\hskip 5mm \rule{2.5mm}{2.5mm}}

\newcommand{\C}{{\mathbb C}}
\newcommand{\N}{{\mathbb N}}
\newcommand{\Z}{{\mathbb Z}}

\newcommand{\proof}{{\em Proof:\ }}

\begin{document}
\newtheorem{thm}{Theorem}[section]
\newtheorem{defs}[thm]{Definition}
\newtheorem{lem}[thm]{Lemma}
\newtheorem{note}[thm]{Note}
\newtheorem{cor}[thm]{Corollary}
\newtheorem{prop}[thm]{Proposition}
\renewcommand{\theequation}{\arabic{section}.\arabic{equation}}
\newcommand{\newsection}[1]{\setcounter{equation}{0} \section{#1}}
\renewcommand{\baselinestretch}{1}

\title{Indefinite boundary value problems on graphs
      }
\author{Sonja Currie \footnote{ Supported by NRF Thutuka grant no. TTK2007040500005} \\ $< Sonja.Currie@wits.ac.za >$ \\ \\
Bruce A. Watson \footnote{ Supported in part by the Centre for Applicable Analysis and
Number Theory. Supported by NRF grant no. FA2007041200006} \\ $<b.alastair.watson@gmail.com>$\\ \\
School of Mathematics\\
University of the Witwatersrand\\
Private Bag 3,\\
P O WITS 2050,\\
South Africa }
\maketitle
\abstract{
 We consider the spectral structure of indefinite second
 order boundary-value
 problems on graphs. A variational formulation for such boundary-value
 problems on graphs is given and we obtain both full and half-range
 completeness results. This leads to a max-min principle and as a
 consequence we can formulate an analogue of Dirichlet-Neumann bracketing
 and this in turn gives rise to asymptotic approximations for the eigenvalues. 
}

\parindent=0in
\parskip=.15in
\newsection{Introduction}
Let $G$ be an oriented graph with finitely many edges, say $K$, each of unit length, having the path-length metric. Suppose that $n$ of the edges have positive weight, $1$, and $K-n$ of the edges have negative weight, $-1$.
We consider the second-order differential equation
\begin{equation}\label{diff}
ly:=-\frac{d^2y}{dx^2} + q(x)y = \lambda By, 
\end{equation}
on $G$, where $q$ is real valued and essentially bounded on $G$ and $By(x) = b(x)y(x)$ with
 \[ b(x):= \left\{\begin{array}{rl}
 1, & \makebox{for $x$ on edges with positive weight.} \\
-1, & \makebox{for $x$ on edges with negative weight.} 
\end{array}\right.\]

At the vertices or nodes of $G$ we impose
formally self-adjoint boundary conditions, see \cite{carlson} for more
details regarding the self-adjointness of boundary conditions.

A variational formulation for a class of indefinite self-adjoint boundary-value problems on graphs is given, see \cite{beals} and \cite{curgus} for background on Sturm-Liouville problems with indefinite weight, and \cite{binding} concerning variational principles in Krein spaces. We then study the nature of the spectrum of this variational problem and obtain both full and half-range completeness results. A max-min principle for indefinite Sturm-Liouville boundary-value problems on directed graphs is then proved which enables us to develop an analogue of Dirchlet-Neumann bracketing for the eigenvalues of the boundary-value problem and consequently to obtain eigenvalue asymptotics.

In parallel to the variational aspects of
boundary-value problems on graphs studied here and on trees
in \cite{solomyak}, the work of Pokornyi and Pryadiev, and
Pokornyi, Pryadiev and Al-Obeid, in \cite{p-p} and \cite{p-p-a}, should be
noted for the extension of Sturmian oscillation theory to
second order operators on graphs.
The idea of approximating the behaviour of eigenfunctions and
eigenvalues for a boundary-value problem
on a graph by the behaviour of associated problems on the individual edges, used here, was studied in the definite case in \cite{ali}, \cite{scbw1} and \cite{vB}.

An extensive survey of the physical systems giving rise to boundary-value problems on graphs
can be found in \cite{kuchment} and the bibliography thereof.
Second order boundary-value problems on finite graphs arise naturally in quantum
mechanics and circuit theory, \cite{Avron, Ger}.
Multi-point boundary-value problems and periodic boundary-value problems
can be considered as particular cases of boundary-value problems
on graphs, \cite{codding}.

In Section \ref{sec2}, the boundary-value problem,
which forms the topic of this paper, is stated and
allowable boundary conditions discussed.
An operator formulation is given along with definitions of
the various function spaces used. A variational reformulation of the boundary-value problem together with the definition of co-normal (elliptic) boundary conditions 
is given in Section \ref{minimax}. Here we also show that a function is a variational eigenfunction if and only if it is a
classical eigenfunction. In Section \ref{nature}, we study the spectrum of the variational problem. The main result of this section is that an eigenfunction is in the positive cone, with respect to the $B$ (indefinite inner product), if and only if the corresponding eigenvalue is positive and similary for the negative cone.
Following the approach used by Beals in \cite{beals} we prove both full and half-range completeness in Section \ref{fhc}, see Theorem \ref{frc} and Theorem \ref{halfrange}. In Section \ref{6}, a max-min characterization of the eigenvalues
of the boundary value problem is given which is then used in Section \ref{7} to obtain a variant of
Dirichlet-Neumann bracketing of the eigenvalues. Hence eigenvalue asymptotics are found. Dirichlet-Neumann bracketing for elliptic partial differential equations can be found in \cite{Cour}.

\section{Preliminaries} \label{sec2}

Denote the edges of the graph G by $e_i$ for $i=1, \dots,K$. As $e_i$ has length $1$, $e_i$ can be considered as the interval $[0,1]$, where $0$ is identified with the initial point of
$e_i$ and $1$ with the terminal point.

We recall, from \cite{scbw1}, the following classes of function spaces:
\begin{eqnarray*}
 \mathcal{L}^2(G)&:=& \bigoplus_{i=1}^K  \mathcal{L}^2(0,1),\\
 \mathcal{H}^m(G)&:=& \bigoplus_{i=1}^K  \mathcal{H}^m(0,1),\quad m=0,1,2,\dots ,\\
 \mathcal{H}^m_o(G)&:=& \bigoplus_{i=1}^K  \mathcal{H}^m_o(0,1),\quad m=0,1,2,\dots ,\\
 \mathcal{C}^\omega(G)&:=& \bigoplus_{i=1}^K  \mathcal{C}^\omega(0,1),\quad \omega=\infty,0,1,2,\dots ,\\
 \mathcal{C}^\omega_o(G)&:=& \bigoplus_{i=1}^K  \mathcal{C}^\omega_o(0,1),\quad \omega=\infty,0,1,2,\dots \ .
\end{eqnarray*}
The inner product on $\mathcal{H}^m(G)$ and $\mathcal{H}^m_0(G)$, denoted $(\cdot,\cdot)_m$,
is defined by
\begin{equation}\label{inner}
(f,g)_m
:=\sum_{i=1}^K \sum_{j=0}^m \int_0^{1} f|_{e_i}^{(j)}\ \bar{g}|_{e_i}^{(j)}\ dt
=:\sum_{j=0}^m \int_G f^{(j)}\ \overline{g}^{(j)}\ dt.
\end{equation}
Note that
$\mathcal{L}^2(G)=\mathcal{H}^0(G)=\mathcal{H}^0_0(G)$. 
For brevity we will write $(\cdot,\cdot)=(\cdot,\cdot)_0$, $\|f\|_m^2=(f,f)_m$ and $\|f\|=\|f\|_0$.

The differential equation (\ref{diff}) on the graph $G$ can be
considered as the system of equations
\begin{equation}\label{diff-1}
-\frac{d^2y_i}{dx^2} + q_i(x)y_i = \lambda b_i(x) y_i,\quad x\in [0,1],\; i=1,\dots ,K,
\end{equation}
where $q_i$, $b_i$ and $y_i$ denote $q|_{e_i}$, $b|_{e_i}$ and $y|_{e_i}$.

As in \cite{scbw1}, the boundary conditions at the node $\nu$ are specified in terms of
the values of $y$ and $y'$ at $\nu$ on each of the incident edges. In particular,
if the edges which start at $\nu$ are $e_i, i\in \Lambda_s(\nu)$,
and the edges which end at $\nu$ are $e_i, i\in \Lambda_e(\nu)$,
then the boundary conditions at $\nu$ can be expressed as
\begin{equation}\label{bc-1}
\sum_{j\in\Lambda_s(\nu)} \left[ \alpha_{ij}y_j+\beta_{ij}{y'}_j\right](0)
+\sum_{j\in\Lambda_e(\nu)} \left[ \gamma_{ij}y_j+\delta_{ij}{y'}_j\right](1)=0,
\quad i=1,\dots ,N(\nu),
\end{equation}
where $N(\nu)$ is the number of linearly independent boundary conditions at node $\nu$.
For formally self-adjoint boundary conditions $N(\nu)=\sharp(\Lambda_s(\nu))+\sharp(\Lambda_e(\nu))$
and $\sum_\nu N(\nu)= 2K$, see \cite{carlson, Naimark} for more details.

Let $\alpha_{ij}=0=\beta_{ij}$ for $i=1,\dots ,N(\nu)$ and
$j\not\in \Lambda_s(\nu)$ and similarly let
$\gamma_{ij}=0=\delta_{ij}$ for $i=1,\dots ,N(\nu)$ and
$j\not\in \Lambda_e(\nu)$. The boundary conditions (\ref{bc-1}) considered over all nodes $\nu$, after possible relabelling,
may thus be written as
\begin{eqnarray}
 \sum_{j=1}^K [\alpha _{ij}y_j(0) + \gamma _{ij}y_j(1)]
&=&0,\quad i=1,\dots,J,\label{gjbc1}\\
 \sum_{j=1}^K [\alpha _{ij}y_j(0)+\beta _{ij}y_j'(0)
+\gamma _{ij}y_j(1)+\delta _{ij}y_j'(1)]&=&0,\quad i=J+1,\dots,2K,
 \label{gjbc2}
\end{eqnarray}
where all possible Dirichlet-like terms are in (\ref{gjbc1}),
i.e. if (\ref{gjbc2}) is written in matrix form then Gauss-Jordan reduction
will not allow any pure Dirichlet conditions linearly independent of (\ref{gjbc1}) to be extracted.

The boundary-value problem (\ref{diff-1})-(\ref{bc-1}) on $G$ can be
formulated as an operator eigenvalue problem in $\mathcal{L}^2(G)$, \cite{agmon, carlson, sho},
for the closed densely defined operator $BL$, where
\begin{equation}
 Lf :=-f''+qf\label{the-operator}
\end{equation}
with domain
\begin{equation}
 \mathcal{D}(L) = \{ f\ |\ f,f' \in AC, Lf\in \mathcal{L}^2(G),
\ f \mbox{ obeying (\ref{bc-1})}\ \}.\label{the-domain}
\end{equation}
The formal self-adjointness of (\ref{bc-1}) relative to $L$ ensures that
$L$ is a closed densely defined self-adjoint operator in $\mathcal{L}^2(G)$, see \cite{H-P, Naimark, Weid}, and that $BL$ is self-adjoint in $H_K$ where $H_K$ is ${\cal L}^2(G)$ with indefinite inner product $[f,g] = (Bf,g)$.

From \cite{scbw1} we have that the operator $L$ is lower semibounded in ${\cal L}^2(G)$.

\section{Variational Formulation}\label{minimax}

In this section we give a, variational formulation for the boundary-value problem (\ref{diff-1})-(\ref{bc-1}) or equivalently
for the eigenvalue problem associated with the operator $BL$. 

\begin{defs} 
(a) Let 
$\mathcal{D}(F) = \{y\in \mathcal{H}^1(G)\ |\ y \makebox{ obeys (\ref{gjbc1})}\}$,
where
$$
\int_{\partial G} y \,d\sigma:= \sum_{i=1}^K [y_i(1)-y_i(0)]
= \int_G y'\,dt.
$$

(b) We say that the boundary conditions on a graph are co-normal or elliptic with respect
to $l$ if  there exists $f$ defined on $\partial G$, such that
$x\in\mathcal{D}(F)$ has
 \[
\int_{\partial G} (fx + x')\overline{y}\, d\sigma = 0, \quad \mbox{for all}\quad
 y \in \mathcal{D}(F)
\]
if and only if $x$ obeys (\ref{gjbc2}).

(c) If the boundary conditions are co-normal and $f$ is as in (b) and ${\cal D}(F)$ is as in (a), then we define the sesquilinear form $F(x,y)$ for $x,y \in {\cal D}(F)$ by 
\begin{equation}\label{varform}
F(x,y):= \int_{\partial G} fx\overline{y}\,d\sigma  + \int_G (x'\overline{y}'
+ xq\overline{y})\,dt.
\end{equation}
\end{defs}

We note that `Kirchhoff', Dirichlet, Neumann and periodic boundary conditions are all co-normal, but this class does not
include all self-adjoint boundary-value problems on graphs. 

The following lemma shows that a function is a variational eigenfunction if and only if it is a
classical eigenfunction.

\begin{lem} \label{regu}
Suppose that (\ref{gjbc1})-(\ref{gjbc2}) are co-normal boundary conditions
with respect to  $l$ of (\ref{diff}). Then  $u \in \mathcal{D}(F)$ satisfies
$F(u,v) = \lambda (Bu,v)$ for all $v \in \mathcal{D}(F)$ if and only if
 $u\in \mathcal{H}^2(G)$ and $u$ obeys (\ref{diff}), (\ref{gjbc1})-(\ref{gjbc2}).
\end{lem}

\proof
 Assume that $u\in \mathcal{H}^2(G)$ and $u$ obeys (\ref{diff}), (\ref{gjbc1})-(\ref{gjbc2}).
 Then for each $v\in \mathcal{D}(F)$
 \begin{eqnarray*}
  F(u,v) &=& \int_{\partial G} fu\overline{v}\, d\sigma + \int_G (u'\overline{v}' + qu\overline{v})\,dt\\
         &=& \int_{\partial G} fu\overline{v}\, d\sigma
                   + \int_G ((u'\overline{v})' - u''\overline{v}+qu\overline{v})\,dt\\
         &=& \int_{\partial G} fu\overline{v}\, d\sigma
                   + \int_G (u'\overline{v})'\,dt +\lambda (Bu,v)\\
         &=& \int_{\partial G} (fu+u')\overline{v}\, d\sigma +\lambda (Bu,v).
 \end{eqnarray*}
 The assumption that (\ref{gjbc1})-(\ref{gjbc2}) are co-normal boundary
conditions with respect to
 $l$ gives that $u\in \mathcal{D}(F)$ and
 $$
\int_{\partial G} (fu+u')\overline{v}\, d\sigma =0,\quad\mbox{for all }
 v\in \mathcal{D}(F),
$$
 completing the proof this in case.

 Now assume $u \in \mathcal{D}(F)$ satisfies $F(u,v) = \lambda (Bu,v)$ for all $v \in \mathcal{D}(F)$.
 As $\mathcal{C}_o^\infty(G)\subset \mathcal{D}(F)$, it follows that
 $$
F(u,v) = \lambda (Bu,v),\quad \mbox{for all } v \in \mathcal{C}_0^{\infty}(G).
$$
 Hence $F(u,\cdot)$ can be extended to a continuous linear functional on
$\mathcal{L}^2(G)$.
 In particular, since $q \in {\cal L}^{\infty}(G)$, this gives that
 $$
\partial u'\in\mathcal{L}^2(G)\subset \mathcal{L}^1_{\rm loc}(G)
$$
 where $\partial$ denotes the distributional derivative.
 Then, by \cite[Theorem 1.6, page 44]{sho},
 $u'\in AC$ and $u''\in \mathcal{L}^1_{\rm loc}(G)$ allowing integration by parts.
 Thus
$$
lu=-u''+qu \in \mathcal{L}^1_{\rm loc}(G)
$$
and consequently $lu=\lambda Bu\in\mathcal{L}^2(G)$.
 Now $q\in \mathcal{L}^\infty(G)$ and
 $\mathcal{D}(F)\subset \mathcal{L}^2(G)$, giving  $u, u''\in \mathcal{L}^2(G)$
 and hence $u\in \mathcal{H}^2(G)$.

 The definition of $\mathcal{D}(F)$ ensures that (\ref{gjbc1}) holds.
 Integration by parts gives
 $$
\int_{\partial G} (fu+u')\bar{y}\ d\sigma =0,\quad \mbox{for all }
 y\in\mathcal{D}(F),
$$
 which, from the definition of $f$ and the constraints on the class of
 boundary conditions,
 is equivalent to $u$ obeying (\ref{gjbc2}). \qed
 
\section{Nature of the spectrum} \label{nature}

The operator $L$ is self-adjoint in ${\cal L}^2(G)$ with spectrum consisting of pure point spectrum and accumulating only at $+ \infty$. In addition, we assume that $L$ is positive definite, thus the spectrum of $L$ may be denoted $0<\rho _1\leq \rho _2\leq \dots$ where $\lim _{n\rightarrow \infty} \rho _n= \infty$. Since $L$ is positive definite and the spectrum consists only of point spectrum, $L^{-1}$ exists and is a compact operator see, \cite[p.24]{sc-thesis}, moreover
\begin{equation}\label{star1}
L^{-1}y(t) = \int_G g(t,\tau )y(\tau )\,d\tau ,
\end{equation}
where $g(t, \tau )$ is the Green's function of $L$. Thus $L^{-1}B$ is a compact operator. Consider the eigenvalue problem
\[\mu y =L^{-1}By, \quad y\in {\cal L}^2(G),\]
where $\mu =\frac{1}{\lambda }$. Since $L^{-1}B$ is compact it has only discrete spectrum except possibly at $\mu =0$ and the only possible accumulation point is $\mu =0$. In addition, $\mu =0$ is not an eigenvalue of $L^{-1}B$ since $0$ is not an eigenvalue of $L^{-1}$. Thus $L^{-1}B$ has countably infinitely many eigenvalues, all non-zero, but accumulating at $0$. From (\ref{star1}) it follows that
\[L^{-1}By(t) = \int_G g(t,\tau )By(\tau )\,d\tau = \int_G \tilde{g}(t,\tau )y(\tau )\,d\tau ,\]
where $\tilde{g}(t,\tau )= g(t,\tau )b(\tau )$. Hence $BL$ has discrete spectrum only, with possible accumluation point at $\infty$ in the complex plane. The spectrum is also countably infinite and, as $0$ is not an eigenvalue of $L$, $0$ is also not an eigenvalue of $BL$.  

\begin{lem}\label{hilb}
The space ${\cal D}(F)$ is a Hilbert space with inner product $F$. The norm generated by $F$ on ${\cal D}(F)$ is equivalent to the $H^1(G)$ norm, making ${\cal D}(F)$ a closed subspace of $H^1(G)$.
\end{lem}

\proof
By (\ref{varform}), \cite[Preliminaries]{scbw1} and the trace theorem, see \cite[p. 38]{agmon} we have that 
there exist constants $K, c>1$ such that
\begin{equation}\label{star2}
\frac{1}{c}||x||^2_{H^1(G)} \leq F(x,x) + K||x||^2 \leq c||x||^2_{H^1(G)}.
\end{equation}
Thus the sesquilinear form $F(x,y) + K(x,y)$ is an inner product on ${\cal D}(F)$. 
From (\ref{star2}) we get directly that
\[\frac{1}{c}(F(x,x) + K||x||^2)\leq ||x||^2_{H^1(G)} \leq c(F(x,x) + K||x||^2),\]
making $F(x,y) + K(x,y)$ and $(x,y)_{H^1(G)}$ equivalent inner products on ${\cal D}(F)$.

We now show that $F(x, y)$ is an inner product on ${\cal D}(F)$ and is equivalent to the inner product $F(x,y) + K(x,y)$ on ${\cal D}(F)$. As $\rho _1$ is the least eigenvalue of $L$ on ${\cal L}^2(G)$, 
\[(Ly,y) \geq \rho _1(y,y) = \rho _1||y||^2,\]
for all $y\in {\cal D}(L)\subset {\cal D}(F)$. Since $F(y,y) = (Ly,y)$, for all $y \in {\cal D}(L)$, we get
\[F(y,y) \geq \rho _1 ||y||^2,\]
for $y \in {\cal D}(L)$.

Now, ${\cal D}(L)$ is dense in ${\cal D}(F)$ for ${\cal D}(F)$ with norm $|||x|||^2 := F(x,x) + K(x,x)$. Thus, by continuity,
\[||y||^2_F := F(y,y) \geq \rho _1 ||y||^2,\]
for all $y \in {\cal D}(F)$, showing that $||\cdot ||_F$ is a norm on ${\cal D}(F)$
and that $F(x, y)$ is an inner product on ${\cal D}(F)$. In addition 
\[\left(1 + \frac{K}{\rho _1}\right)||y||_F^2 = F(y,y) + \frac{K}{\rho _1} F(y,y) \geq F(y,y) + K(y,y)\geq F(y,y)= ||y||_F^2,\]
where $K$ is as given above.
Thus $F(x,y) + K(x,y)$ and $F(x,y)$ are equivalent inner products on ${\cal D}(F)$ and since $F(x,y) + K(x,y)$ and $(x,y)_{H^1(G)}$ are equivalent inner products on ${\cal D}(F)$ we have that $F(x,y)$ and $(x,y)_{H^1(G)}$ are equivalent inner products on ${\cal D}(F)$.

We now show that, with the $F$ inner product, ${\cal D}(F)$ is a Hilbert space. For this, we need only show that ${\cal D}(F)$ is closed in $H^1(G)$. The map $\hat{T}: H^1(G) \rightarrow \C^J$ given by $$\hat{T} : y \rightarrow  \left(\sum_{j=1}^K [\alpha _{ij}y_j(0) + \gamma _{ij}y_j(1)]\right)_{i= 1, \dots, J},$$ is continuous by the trace theorem, see \cite{agmon}, and thus the kernel of $\hat{T}$, Ker$(\hat{T}) = {\cal D}(F)$ is closed. \qed   


\begin{thm}\label{realspec}
The spectrum of (\ref{diff}), (\ref{gjbc1})-(\ref{gjbc2}) is real and all eigenvalues are semi-simple.
\end{thm}

\proof
As ${\cal D}(L)$ is dense in ${\cal D}(F)$, $L$ is a densely defined operator in ${\cal D}(F)$.
Now $F(x,y) := (Lx,y)$ for all $x \in {\cal D}(L)$ and $y \in {\cal D}(F)$. 

Let $\tilde{L}:= L^{-1}B$, then $\tilde{L}: {\cal L}^2(G) \rightarrow {\cal D}(L)$ and is thus a map from ${\cal D}(F)$ to ${\cal D}(L)$.

Since $B$ and $L$ are self adjoint in $ \mathcal{L}^2(G)$ we get
\begin{eqnarray*}
F(\tilde{L}x,y) &=& F(L^{-1}Bx,y) \\
&=& (Bx,y) \\
&=& (x,By)\\
&=& \overline{(By,x)}\\
&=& \overline{F(\tilde{L}y, x)}\\
&=& F(x, \tilde{L}y).
\end{eqnarray*}
for $x,y \in \mathcal{D}(F)$.

So $\tilde{L}$ is self adjoint in $\mathcal{D}(F)$ (with respect to $F$). Thus, in ${\cal D}(F)$, $\tilde{L}$ has only real spectrum  and all eigenvalues are semi-simple. Therefore, by Lemma \ref{regu}, the pencil $Lx=\lambda Bx$ has only real spectrum and all eigenvalues are semi-simple. \qed

Let
\begin{equation}\label{Kinner}
[f,g]
:=\sum_{i=1}^n \int_0^{1} f|_{e_i}\ \bar{g}|_{e_i}\, dt - \sum_{i=n+1}^K \int_0^{1} f|_{e_i}\ \bar{g}|_{e_i}\, dt = (Bf,g), 
\end{equation}
then $ \mathcal{L}^2(G)$, with the indefinite inner product given by (\ref{Kinner}), is a Krein space which we denote by $H_K$.

We now define the positive, $C^+$, and negative, $C^-$, cones of $H_K$ by
\[C^+ := \{y \in H_K\,|\,[y,y]>0\},\]
\[C^- := \{y \in H_K\,|\,[y,y]<0\}.\]

\begin{thm}
For $L$ positive definite in ${\cal L}^2(G)$ and $y$ an eigenfunction of (\ref{diff}), (\ref{gjbc1})-(\ref{gjbc2}) corresponding to the eigenvalue $\lambda$ we have
$y\in C^+$ if and only if $\lambda >0$, and $y\in C^-$ if and only if $\lambda <0$.
\end{thm}

\proof
Let $y$ be an eigenfunction corresponding to $\lambda $.
Using the fact that any element, $y$, of $H_K$ may be written in the form $y = \{f,g\}$ or $y =f \oplus g$, where $f = (y|_{e_1}, \dots, y|_{e_n})$ has $n$ components and $g = (y|_{e_{n+1}}, \dots, y|_{e_K})$ has $K-n$ components, we get that
\[C^+ = \{ \{f,g\} \,|\,||f||^2_{{\cal L}^2(G^+)}> ||g||^2_{{\cal L}^2(G^-)}\},\] and
\[C^- = \{ \{f,g\} \,|\,||f||^2_{{\cal L}^2(G^+)}< ||g||^2_{{\cal L}^2(G^-)}\}.\]
Here $G^+$ denotes the subgraph of $G$ where the weights are positive and $G^-$ denotes the subgraph of $G$ where the weights are negative.

Since $L>0$ and $y = \{f,g\}$,
\[0<(Ly,y) = (\lambda By,y) = \lambda [y,y] = \lambda (||f||^2_{{\cal L}^2(G^+)} - ||g||^2_{{\cal L}^2(G^-)}).\]
Hence, $y\in C^+$ if and only if $\lambda > 0$, and $y\in C^-$ if and only if $\lambda < 0$. \qed

\section{Full and half-range completeness}\label{fhc}



In this section we prove both half and full range completeness of the eigenfunctions of (\ref{diff}), (\ref{gjbc1})-(\ref{gjbc2}). In the case presented here the proof is simpler than that of Beals \cite{beals}, but it is assumed that the problem is left definite, i.e. $L$ is a positive operator.


Recall that, by Lemma \ref{hilb}, ${\cal D}(F)$ is a Hilbert space. Define $$\tilde{F}[u](v):= F(u,v)$$ then $\tilde{F} : {\cal D}(F) \longrightarrow {\cal D}(F)'$, where ${\cal D}(F)'$ is the conjugate dual of ${\cal D}(F)$, i.e. the space of continuous conjugate-linear maps from ${\cal D}(F)$ to $\C$. 

\begin{lem}
$\tilde{F}$ is an isomorphism from ${\cal D}(F)$ to ${\cal D}(F)'$.
\end{lem}

\proof
If $F(u_1, v)=F(u_2, v)$, for all $v\in {\cal D}(F)$,  then $u_1=u_2$ since $F$ is an inner product on ${\cal D}(F)$, see Lemma \ref{hilb}. Thus $\tilde{F}$ is one to one. 

Now, for ${\hat v} \in {\cal D}(F)'$ we have that ${\hat v}(x) = F(v,x)$ for some $v \in {\cal D}(F)$ by the Theorem of Riesz, \cite{Rudin}. So ${\hat v}(x) = \tilde{F}[v](x)$ giving that $\tilde{F}[v] = {\hat v}$. Hence $\tilde{F}$ is onto. 

Also $\tilde{F}$ and $\tilde{F}^{-1}$ are everywhere defined maps on a Hilbert space and are thus continuous as a consequence of the principle of uniform boundedness (Banach Steinhaus theorem), \cite{Rudin}.

So $\tilde{F}$ is an isomorphism from ${\cal D}(F)$ to ${\cal D}(F)'$. \qed

Define $T[u](v) := (Bu,v)$ for $u,v \in {\cal D}(F)$. Then $T: {\cal D}(F) \longrightarrow {\cal D}(F)'$ is compact since ${\cal D}(F)$ is compactly embedded in ${\cal L}^2(G)$ and $Bu \in {\cal L}^2(G)$ with the mapping $Bu \mapsto (Bu, \cdot )$ from ${\cal L}^2(G)$ to ${\cal L}^2(G)'$ continuous. Thus $S:= \tilde{F}^{-1}T$ is a compact map with $S : {\cal D}(F) \longrightarrow {\cal D}(F)$.

\begin{lem}
The compact operator $S$ on ${\cal D}(F)$ is self-adjoint with respect to the inner product $F$.
\end{lem} 

\proof
For $u,v \in {\cal D}(F)$
\[F(Su,v) = \tilde{F}[Su](v)=T[u](v) = (Bu,v) = (u,Bv).\]
Similarly
   \[\overline{(Bv,u)}= \overline{F(Sv,u)} = F(u,Sv). \qed\]
  
As $S$ is a compact self-adjoint operator on ${\cal D}(F)$ and as $0$ is not an eigenvalue of $S$, the eigenfunctions, $(u_n)$, of $S$, with eigenvalues $(\lambda _n^{-1})$, can be chosen so that $(u_n)$ is an orthonormal basis for ${\cal D}(F)$. 

{\bf Note:} The equation $Su_n = \lambda _n^{-1}u_n$ is equivalent to the equation $Lu_n = \lambda _n Bu_n$, in the sense that if
\[\lambda _n Su_n = u_n,\] then, by the definition of $S$,
\[\lambda _n(\tilde{F}^{-1}T)u_n = u_n.\] Applying $\tilde{F}$ to the above gives 
\[\lambda _nTu_n = \tilde{F}u_n.\] Thus
\[\lambda _nT[v](u_n) = \tilde{F}[v](u_n),\] for all $v\in {\cal D}(F)$. From the definition of $T$, this gives
\[\lambda _n (Bv, u_n) = \tilde{F}[v](u_n).\] Hence
\[\lambda _n (Bv, u_n) = F(v,u_n)\] for all $v\in {\cal D}(F)$. Using Lemma \ref{regu} we we obtain that
\[\lambda _n (Bv, u_n) = (v,Lu_n).\] Therefore
\[(v,\lambda _n Bu_n -Lu_n) =0,\] for all $v\in {\cal D}(F)$, and by the density of ${\cal D}(F)$ in ${\cal L}^2(G)$, this yields
\[Lu_n = \lambda _n Bu_n.\]
It is easy to show that if $Lu_n = \lambda _n Bu_n$, then $Su_n = \lambda _n^{-1}u_n$.

In summary, we have the following theorem:

\begin{thm}[Full range completeness]\label{frc}
The eigenfunctions $(y_n)$ of (\ref{diff}), (\ref{gjbc1})-(\ref{gjbc2}) form a Riesz basis for ${\cal L}^2(G)$ and can be chosen to form an orthonormal basis for ${\cal D}(F)$ (with respect to the $F$ inner product). In addition $(y_n)$ is orthogonal with respect to $[\cdot ,\cdot ]$.
\end{thm}

\proof
Since $S$ is a compact self-adjoint operator on the Hilbert space ${\cal D}(F)$, the eigenvectors can be chosen to form an orthonormal basis in ${\cal D}(F)$. As shown in the note above the variational eigenfunctions coincide with those of $L^{-1}B$ (with eigenvalues mapped by $\lambda \mapsto \frac{1}{\lambda }$ and where $0$ is not in the point spectrum).Thus the eigenfunctions of $L^{-1}B$ can be chosen to form an orthonormal basis for ${\cal D}(F)$ and as ${\cal D}(F)$ is dense in ${\cal L}^2(G)$ they form a Riesz basis for ${\cal L}^2(G)$.

Finally, if $(y_n)$ is an orthonormal basis of ${\cal D}(F)$ of eigenfunctions then
\[\delta _{n,m} = F(y_n,y_m) = (\lambda _nBy_n, y_m) = \lambda _n(By_n,y_m) = \lambda _n[y_n,y_m].\]
Hence $(y_n)$ is orthogonal with respect to $[\cdot ,\cdot ]$. \qed

Let $P_{\pm }$ be the positive and negative spectral projections of $S$. Note that Ker$(S) = \{0\}$. The projections, $P_{\pm }$, are then defined by the property
\[P_{\pm }u_n = \left\{ \begin{array}{cl}
u_n, & \pm \lambda _n>0\\
0, & \pm \lambda _n<0
\end{array}\right.,\]
hence
\[|S| = S(P_+ - P_-) = (P_+ - P_-)S.\]
On ${\cal D}(F)$ we introduce the inner product $(u,v)_{S} = F(|S|u,v)$ with related norm $||u||_S = (u,u)_S^{\frac{1}{2}}$.

We must now show that this norm is equivalent to the ${\cal L}^2(G)$ norm, $||u||=(u,u)^{\frac{1}{2}}$. 

The operator $B$ is a self-adjoint operator in ${\cal L}^2(G)$ and $B$ has spectral projections $Q_{\pm }$, where
\[Q_{\pm }u(x) = \left\{ \begin{array}{cl}
u(x), & b(x)=\pm 1\\
0, & b(x) = \mp 1
\end{array}\right..\]

Thus $|B| = I = B(Q_+ + Q_-) = (Q_+ + Q_-)B$ is just the identity map, and $|T|$ is the map from ${\cal D}(F)$ to ${\cal D}(F)'$ induced by $|B|$, i.e. $|T|[u](v) = (u,v)$.  
But $T[u](v) := (Bu,v)$ for all $u,v \in {\cal D}(F)$, and thus can be extended to $u,v \in {\cal L}^2(G)$, i.e.
\[T: {\cal L}^2(G) \rightarrow {\cal L}^2(G)'\hookrightarrow {\cal D}(F)'.\]
In this sense $TQ_{\pm} : {\cal L}^2(G) \rightarrow {\cal D}(F)'$ is compact.

Also $T(Q_+ + Q_-)[u](v) = (B(Q_+ +Q_-)u,v) = (u,v) = |T|[u](v)$ for all $u,v \in {\cal L}^2(G)$ and thus for $u,v\in {\cal D}(F)$. We now observe that $Q_{\pm }'T : {\cal D}(F)\rightarrow {\cal D}(F)'$, using the extension of $T$ to ${\cal L}^2(G)$, is well defined as $Q_{\pm }'T[u](v) = T[u](Q_{\pm }v) = (Bu,Q_{\pm }v) = (Q_{\pm }Bu,v) = (BQ_{\pm }u,v)$ making  $TQ_{\pm } = Q_{\pm }'T$. Hence
\[|T| = T(Q_+ - Q_-) = (Q_+' -Q_-')T.\]

\begin{thm}\label{eqivnorms}
The norms $||\cdot ||_S$ and $||\cdot ||$ are equivalent on ${\cal D}(F)$.
\end{thm}

\proof
Considered as an operator in the subspace $P_+({\cal D}(F))$, $S$ is a positive operator. Let $y\in {\cal D}(L)$. Since $L$ is a positive operator and ${\cal D}(F)$ is compactly embedded in ${\cal L}^2(G)$ we have that there exists some constant $C>0$ such that
\begin{equation}\label{p}
(Ly,y)=F(y,y)\geq C(y,y),
\end{equation}
for all $y \in {\cal D}(L)$. Also
\begin{equation}\label{q}
||Q_+y||^2 \leq ||y||^2.
\end{equation}
Combining (\ref{p}) and (\ref{q}) we obtain that
\begin{equation}\label{c}
C||Q_+y||^2 \leq  C(y,y) \leq (Ly,y),
\end{equation}
for $y \in {\cal D}(L)$. Let $(y_n)$ be an orthonormal basis of eigenfunctions of $S$ in ${\cal D}(F)$ where $y_n$ has eigenvalue $\lambda _n$ with $0<\lambda _1<\lambda _2< \dots$ and $0>\lambda _{-1}>\lambda _{-2}> \dots.$ Now
\[P_+({\cal D}(F)) = \overline{<y_1, y_2, \dots >},\]
and $Ly_n = \lambda _nBy_n$ for all $n=1,2, \dots$. 

Let $y\in P_+({\cal D}(L))$ then $y = \sum_{n=1}^{\infty} \alpha _ny_n$ where $\alpha _n \in \C, n\in \N$. From (\ref{c}) we have that 
\[||Q_+y||^2 \leq \frac{1}{C}(Ly,y).\]
Using the orthogonality of $(y_n)$ we get
\[ \frac{1}{C}(Ly,y) = \sum_{n=1}^{\infty} |\alpha _n|^2 \frac{\lambda _n}{C}(By_n, y_n),\] thus
\[||Q_+ y||^2 \leq  \frac{\lambda _1}{C} \sum_{n=1}^{\infty} |\alpha _n|^2 (By_n, y_n).\]

But 
\[\sum_{n=1}^{\infty} |\alpha _n|^2 (By_n, y_n) = (By,y),\] hence
\begin{eqnarray*}
||Q_+ y||^2 &\leq & \frac{\lambda _1}{C}(By,y)\\
&=& \frac{\lambda _1}{C} T[y](y)\\
&=& \frac{\lambda _1}{C} \tilde{F}[Sy](y)\\
&=& \frac{\lambda _1}{C} F(Sy,y)\\
&=& \frac{\lambda _1}{C} F(|S|y,y).
\end{eqnarray*}
So \[||Q_+ y||^2 \leq \frac{\lambda _1}{C} ||y||_S^2\] and setting $\sqrt{\frac{\lambda _1}{C}} := k >0$ gives
\begin{equation}\label{in1}
||Q_+ y|| \leq k ||y||_S. 
\end{equation}
Similarly 
\[||Q_- y||^2 \leq \frac{\lambda _1}{C} ||y||_S^2\] i.e.
\begin{equation}\label{in2}
||Q_- y|| \leq k ||y||_S. 
\end{equation} 
Since ${\cal D}(L)$ is dense in ${\cal D}(F)$, (\ref{in1}) and (\ref{in2}) hold on all $P_+({\cal D}(F))$, so as $||y||^2 = ||Q_+y||^2 + ||Q_-y||^2$ we have $||y|| \leq \sqrt{2}k ||y||_S$ for all $y \in P_+({\cal D}(F))$.

Working on $P_-({\cal D}(F))$ yields a similar estimate but with $\lambda _1$ replaced by $-\lambda _{-1}$. Thus there exists a constant $C_1>0$ so that for all $y\in {\cal D}(F)$,
\begin{equation}\label{in3}
||y||\leq C_1||y||_S.
\end{equation}

To obtain (\ref{in4}), the reverse of (\ref{in3}), we observe that
\[||y||_S^2 =F(|S|y,y) = F((SP_+ - SP_-)y,y).\]
But $SP_{\pm } = P_{\pm }S$ so
\begin{eqnarray*}
||y||_S^2 &=& F(Sy, P_+y - P_-y)\\
&=& \tilde{F}[Sy](P_+y - P_-y)\\
&=& T[y](P_+y - P_-y)\\
&=& |T|[Q_+y - Q_-y](P_+y - P_-y).
\end{eqnarray*}
Using H\"older's inequality we obtain that
\[|T|[Q_+y - Q_-y](P_+y - P_-y) \leq  ||Q_+y - Q_-y||\,||P_+y - P_-y||.\]
Thus 
\[||y||_S^2 \leq  ||Q_+y - Q_-y||\,||P_+y - P_-y|| = ||y||\,||P_+y - P_-y||.\]
By (\ref{in3})
\[||y||_S^2 \leq C_1||y||\,||P_+y - P_-y||_S.\]
Now
\begin{eqnarray*}
||P_+y - P_-y||_S&=& F(|S|(P_+ - P_-)y, (P_+ - P_-)y)\\
&=& F(Sy, (P_+ - P_-)y)\\
&=& F((P_+-P_-)Sy,y)\\
&=& F(|S|y,y), 
\end{eqnarray*}
giving \[||y||_S^2\leq C_1||y||\,||y||_S,\]
therefore 
\begin{equation}\label{in4}
||y||_S\leq C_1||y||.
\end{equation}
 Combining (\ref{in3}) and (\ref{in4}) gives
 \[\frac{1}{C_1}||y||_S \leq ||y|| \leq C_1||y||_S\] and thus the two norms are equivalent in ${\cal D}(F)$. \qed

Let $H_S$ be the completion of ${\cal D}(F)$ with respect to $||\cdot ||_S$.

\begin{thm}[Half-range completeness]\label{halfrange}
For $Q_+$ and $Q_-$ as previously defined \\ $\{Q_+y_n, \lambda _n>0\}$ is a Riesz basis for ${\cal L}^2(G^+)$ and $\{Q_-y_n, \lambda _n<0\}$ is a Riesz basis ${\cal L}^2(G^-)$.
\end{thm}

\proof
To prove the half-range completeness we show that $\{y_n, \lambda _n>0\}$ and $\{y_n, \lambda _n<0\}$ are Riesz bases for $Q_+P_+(H_S)$ and $Q_-P_-(H_S)$ respectively via showing that $V:= Q_+P_+ + Q_-P_-$ is an isomorphism  from $H_S$ to ${\cal L}^2(G)$, see \cite{beals}.

Let $u,v \in {\cal D}(F)$, then 
\begin{equation}\label{p1}
(Q_{\pm }u, P_{\pm }v)_S = (Q_{\pm }u, P_{\pm }v)
\end{equation} and
\begin{equation}\label{p2}
(Q_{\pm }u, P_{\mp }v)_S = -(Q_{\pm }u, P_{\mp }v).
\end{equation}
To see this, as $S$ is self-adjoint with respect to $F$ so is $|S|$, we have, for example,
\begin{eqnarray*}
(Q_+u, P_-v)_S &=& F(|S|Q_+u, P_-v)\\
&=& F(Q_+u, |S|P_-v)\\
&=& F(Q_+u, S(P_+ - P_-)P_-v)\\
&=& F(SQ_+u, -P_-v)\\
&=& -F(SQ_+u, P_-v)\\
&=& -(Q_+u, P_-v),
\end{eqnarray*}
because $F(SQ_+u, P_-v)= (BQ_+u, P_-v)$ and $Q_+u(x) =0$ when $b(x)=-1$ and $Q_+u(x) = u(x)$ when $b(x)=1$.

Now, as $P_{\pm }$ are self-adjoint with respect to $[\cdot ,\cdot ]$,
\begin{eqnarray*}
||u||_S^2 &=& F(|S|u,u)\\
&=& F((P_+ -P_-)Su,u)\\
&=& F(Su, (P_+ -P_-)u)\\
&=& (Bu, (P_+ -P_-)u)\\
&=& ((Q_+ - Q_-)u, (P_+ -P_-)u)\\
&=& (Q_+u,P_+u) + (Q_-u,P_-u) - (Q_+u,P_-u) - (Q_-u,P_+u).
\end{eqnarray*}
For $u \in {\cal D}(F)$,
\begin{eqnarray*}
||Vu||^2 &=& (Q_+P_+u,Q_+P_+u) + (Q_-P_-u, Q_-P_-u) + (Q_-P_-u,Q_+P_+u) + (Q_+P_+u,Q_-P_-u)\\
&=& (Q_+P_+u,Q_+P_+u) + (Q_-P_-u, Q_-P_-u)\\
&=& (Q_+(I-P_-)u,(I-Q_-)P_+u) + (Q_-(I-P_+)u, (I-Q_+)P_-u)\\
&=& (Q_+u,P_+u) - (Q_+P_-u, P_+u) + (Q_-u,P_-u) - (Q_-P_+u, P_-u)\\
&=& ||u||_S^2  + (Q_+u,P_-u) + (Q_-u,P_+u) - (Q_+P_-u, P_+u)- (Q_-P_+u, P_-u).
\end{eqnarray*}
Setting $W:= Q_+P_- + Q_-P_+$, since $Q_+ -Q_- = B$ and $P_{\pm }$ are self-adjoint and orthogonal with respect to $[\cdot ,\cdot ]$, we obtain
\begin{eqnarray*}
||Vu||^2 &=& ||u||_S^2 +(Q_-P_+u,Q_-P_+u) + (Q_+P_-u,Q_+P_-u)\\
&=& ||u||_S^2 + ||Wu||^2.
\end{eqnarray*}

As $||\cdot ||$ and $||\cdot ||_S$ are equivalent norms on ${\cal D}(F)$, the above equality holds for $u\in H_S$ and shows that the bounded operator $V$ has closed range and kernel $(0)$.

Equations (\ref{p1}) and (\ref{p2}) show that, as mappings from $H_S$ to ${\cal L}^2(G)$, $V$ and $W$ have adjoints
$V^* = P_+Q_+ + P_-Q_-$ and $W^*= -P_+Q_- - P_-Q_+$. But $V^*$ and $W^*$ obey, by the same reasoning as above,
\begin{equation}\label{n2}
||V^*u||_S^2 = ||W^*u||_S^2 + ||u||^2.
\end{equation}
Thus $V^*$ is one to one and therefore $V$ is an isomorphism. Hence we have proved the theorem. \qed

\section{Max-Min Property}\label{6}

In this section we give a maximum-minimum characterization for the
eigenvalues of indefinite boundary-value problems on
graphs. We refer the reader to
\cite[page 406]{Cour} and  \cite{Wein} where analogous results for
partial differential operators were
considered. 

In the following theorem
$\{ v_1,\dots,v_{n}\}^\perp$ will denote the orthogonal complement with respect to $[\cdot ,\cdot ] = (B\cdot ,\cdot )$
of $\{v_1,\dots,v_{n}\}$.
In addition, as is customary, it will be assumed that the eigenvalues,
$\lambda_n>0$, $n\in \N$, of (\ref{diff}), (\ref{gjbc1})-(\ref{gjbc2}), are listed in increasing order and repeated according to
multiplicity, and that the eigenfunctions, $y_n$, are chosen
 so as to form a complete orthonormal family in $\mathcal{L}^2(G)\cap C^+$. More precisely, as in Theorem~\ref{frc}, $(y_n)$, $n\in \Z\setminus \{0\}$ can be chosen so as to form an orthonormal basis for ${\cal D}(F)$ and thus for ${\cal L}^2(G)$ with respect to $B$. In particular $(y_n)_{n\in \N}$ is then an orthonormal basis for ${\cal L}^2(G)\cap C^+$ with respect to $B$ (i.e. $[\cdot ,\cdot ]$). 
The case of ${\cal L}^2(G)\cap C^-$ is similar, so for the remainder of the paper we will restrict ourselves to ${\cal L}^2(G)\cap C^+$.

\begin{thm}\label{mm}
 Suppose $(L\varphi ,\varphi )>0$ for all $\varphi \in {\cal D}(L)\setminus \{0\}$, and for $v_j \in \mathcal{L}^2(G)\cap C^+, j=1,2,\dots$, let
 \begin{equation}\label{dmin}
  d_{n+1}(v_1, \dots, v_{n})
   =\inf\, \left\{\left.\frac{F(\varphi ,\varphi )}{(B\varphi ,\varphi )}\ \right|\
    \varphi \in \{ v_1,\dots,v_{n}\}^\perp \cap D(F)\setminus \{0\}, (B\varphi ,\varphi )>0\right\}.
 \end{equation}
 Then
\begin{equation}\label{A}
\lambda_{n+1} = \sup\,\{d_{n+1}(v_1,\dots ,v_{n})\ |\ v_1,\dots ,v_{n} \in
\mathcal{L}^2(G)\cap C^+\},
\end{equation}
for $n=0,1,\dots, $ and this maximum-minimum is attained if and only if $\varphi=y_{n+1}$ and
$v_i = y_i$, $i=1,\dots ,n$, where $y_j$ is an eigenfunction of $L$ with eigenvalue $\lambda _j$, and $(y_j)$ is a $B$-orthogonal family.
\end{thm}

\proof
 Let $v_1, \dots, v_{n} \in \mathcal{L}^2(G)\cap C^+$.
 As $\text{span}\{y_1, \dots, y_{n+1}\}$ is $n+1$ dimensional and
$\text{span}\{v_1,\dots, v_{n}\}$ is at most $n$
 dimensional there exists
 $\varphi$ in $\text{span}\{y_1, \dots, y_{n+1}\}\setminus \{0\}$
 having
 $$
(B\varphi , v_i) = 0,\quad\mbox{for all}\quad i=1,\dots ,n.
$$
 In particular, this ensures that $\varphi\in\mathcal{D}(F)$ as each
$y_i$ is in $\mathcal{D}(F)$.

 Denote ${\varphi=\sum_{k=1}^{n+1} c_ky_k}$, then
  \begin{eqnarray*}
    F(\varphi,\varphi)&=& \sum_{i,k=1}^{n+1} c_i\bar{c}_k F(y_i,y_k)\\
    &=&\sum_{i=1}^{n+1} |c_i|^2 F(y_i,y_i)\\
    &=& \sum_{i=1}^{n+1} |c_i|^2 (Ly_i,y_i)\\
    &=& \sum_{i=1}^{n+1} |c_i|^2(\lambda _iBy_i,y_i)\\
    &=&\sum_{i=1}^{n+1} |c_i|^2 \lambda _i (By_i,y_i)\\
    & \le & \lambda_{n+1} \sum_{i=1}^{n+1} |c_i|^2 (By_i,y_i)\\
    &=& \lambda_{n+1} (B\varphi ,\varphi ),
 \end{eqnarray*}
 thus showing that
 $$
d_{n+1}(v_1, \dots, v_{n}) \leq \lambda_{n+1} \quad\mbox{for all }
 v_1,\dots,v_{n}\in\mathcal{L}^2(G)\cap C^+.
$$
Hence
 $$\sup\,\{d_{n+1}(v_1,\dots ,v_{n})\,|\,v_1,\dots ,v_{n} \in
\mathcal{L}^2(G)\cap C^+\} \le \lambda_{n+1}.$$

Now suppose $\lambda _{n+1} > d_{n+1}(y_1, \dots, y_n)$. Then there exists $u\in {\cal D}(F)\setminus \{0\}$, $u\in \{y_1, \dots, y_n\}^{\perp }$, such that $B(u,u)=1$ and 
\begin{equation}\label{new1}
F(u,u) < d_{n+1}(y_1, \dots, y_n) + \frac{1}{2}\left(\lambda _{n+1}-d_{n+1}(y_1, \dots, y_n)\right).
\end{equation}

By Theorem \ref{frc} we can write $\displaystyle{u = \sum_{j\not\in \{0, \dots, n\}} \alpha _j y_j}$. Therefore
\begin{eqnarray*}
 F(u,u)&=& \sum_{i,j\not\in \{0, \dots, n\}} \alpha _i\bar{\alpha }_j F(y_i,y_j)\\
    &=&\sum_{i`\not\in \{0, \dots, n\}} |\alpha _i|^2 F(y_i,y_i)\\
    &=& \sum_{i\not\in \{0, \dots, n\}} |\alpha _i|^2 (Ly_i,y_i)\\
    &=& \sum_{i\not\in \{0, \dots, n\}} |\alpha _i|^2(\lambda _iBy_i,y_i).
\end{eqnarray*}
Now as $\lambda _i(By_i,y_i)= F(y_i,y_i)>0$ for all $i$, we have 
\begin{eqnarray*}
 F(u,u) &=&\sum_{i>n} |\alpha _i|^2 \lambda _i (By_i,y_i) + \sum_{i\leq -1} |\alpha _i|^2 \lambda _i (By_i,y_i)\\
    &\geq & \sum_{i>n} |\alpha _i|^2 \lambda _i (By_i,y_i)\\
    &\geq & \lambda _{n+1} \sum_{i>n} |\alpha _i|^2 (By_i,y_i)\\
    &=& \lambda _{n+1}\left(B\sum_{i>n}\alpha _iy_i , \sum_{j>n} \alpha _jy_j \right)\\
    &=& \lambda_{n+1} (BP_+u,P_+u ).
    \end{eqnarray*}
    
Combining the above with (\ref{new1}) and noting that $(Bu,u)=1$, gives
\[\lambda_{n+1} - \frac{1}{2}\left(\lambda _{n+1}-d_{n+1}(y_1, \dots, y_n)\right) >  \lambda_{n+1} (BP_+u,P_+u ).\]
Thus
\[ (Bu,u) - \frac{\lambda _{n+1} -d_{n+1}(y_1, \dots, y_n)}{2\lambda _{n+1}} = 1 - \frac{\lambda _{n+1} -d_{n+1}(y_1, \dots, y_n)}{2\lambda _{n+1}} > (BP_+u,P_+u ).\]
Using the self-adjointness of the projections $P_{\pm }$ with respect to $[\cdot ,\cdot ]$ now gives
\[(BP_-u,P_-u)>\frac{\lambda _{n+1} -d_{n+1}(y_1, \dots, y_n)}{2\lambda _{n+1}} > 0.\]

But $P_-u \in C^-$, so we have a contradiction and therefore $\lambda _{n+1} \leq d_{n+1}(y_1, \dots, y_n)$.

We have shown that $\lambda _{n+1} = d_{n+1}(y_1, \dots, y_n)$, (\ref{A}) holds and $d_{n+1}$ attains its supremum for $(y_1, \dots, y_n)$. Also a direct computation gives $F(y_{n+1}, y_{n+1}) = \lambda _{n+1}(By_{n+1}, y_{n+1})$.

It remains to be shown that if $u \in {\cal D}(F)$ is such that the maximum is attained for $u, v_1, \dots, v_n$ then $u$ is an eigenfunction with eigenvalue $\lambda =d_{n+1}(v_1, \dots, v_n)$.

Let $u\in {\cal D}(F)$ with $(Bu,u) =1$ and 
 $$
J(\varphi,\epsilon)=\frac{F(u+\epsilon\varphi, u+\epsilon\varphi)}{(B(u+ \epsilon \varphi ), u+ \epsilon \varphi )}
\quad\mbox{for all }   \varphi \in \mathcal{D}(F),
\epsilon \in \mathbb{R},\,\, |\epsilon |\,\, \makebox{small}.
$$
 Differentiation with respect to $\epsilon$ of $J(\varphi,\epsilon)$ gives
 $$
0=\frac{\partial}{\partial\epsilon}J(\varphi,\epsilon)|_{\epsilon=0}
   =2\Re[F(\varphi,u)-d_{n+1}(v_1,\dots,v_{n})(B\varphi,u)],
$$
 for all $\varphi \in \mathcal{D}(F)$ and $(Bu,u)=1$.
 Since everything in the above expression is real we obtain that
 \begin{equation}\label{FU}
 F(\varphi,u)=d_{n+1}(v_1,\dots,v_{n})(B\varphi,u),
 \end{equation}
 for all $\varphi \in \mathcal{D}(F)$ and $(Bu,u)=1$.
 
 Now $F(u,u)>0$ therefore $d_{n+1}(v_1,\dots,v_{n})(Bu,u)>0$ which, since $(Bu,u)=1$, gives $d_{n+1}(v_1,\dots,v_{n})>0$. From (\ref{FU}), for $\varphi \in {\cal C}^{\infty}_{0}(G)$,
 we get that 
 \[(L\varphi ,u) - d_{n+1}(v_1,\dots,v_{n})(B\varphi,u) =0, \]
 giving
 \[(\varphi , (l-d_{n+1}(v_1,\dots,v_{n})B)u) =0 .\]
 
Hence, by the proof of Lemma \ref{regu}, $u\in H^2(G)\cap {\cal D}(F)$ and obeys (\ref{diff}) and (\ref{gjbc1}). We must still show that $u$ obeys the boundary condition (\ref{gjbc2}).  
 
 From the proof of Lemma \ref{regu} we see that, for $\varphi \in {\cal D}(F)$,
 \[F(u, \varphi ) = \int_{\partial G}(fu  + u ')\bar{\varphi }d\, \sigma + d_{n+1}(v_1, \dots, v_n)(Bu ,\varphi ). \]
 
 This together with (\ref{FU}) gives that
 \begin{equation}\label{phi0}
 0 = \int_{\partial G}(fu  + u ')\bar{\varphi }d\, \sigma 
 \end{equation}
for all $\varphi \in {\cal D}(F)$. 

As, (\ref{phi0}) holds for all $\varphi \in {\cal D}(F)$, $u$ obeys (\ref{gjbc2}), giving that $u$ is an eigenfunction of (\ref{diff}), (\ref{gjbc1})-(\ref{gjbc2}) with
 eigenvalue $\lambda = d_{n+1}(v_1, \dots, v_n)$. \qed

 \newsection{Eigenvalue Bracketing and Asymptotics}\label{7}

If the boundary conditions (\ref{gjbc1})-(\ref{gjbc2}) are replaced by the Dirichlet condition $y=0$ at each node of $G$,
i.e.
\begin{eqnarray}\label{dir}
 y_i(1)=0 \quad \text{and} \quad y_i(0)=0, \quad i=1,\dots, K,
\end{eqnarray}
then the graph $G$ becomes disconnected with each edge $e_i$ becoming a component sub-graph, $G_i$, with Dirichlet
boundary conditions at its two nodes (ends). The boundary value problem on each sub-graph $G_i$ is equivalent
to a Sturm-Liouville boundary value problem on $[0,1]$ with Dirichlet boundary conditions. Depending on whether the edge has positive or negative weight the resulting boundary value problem is
\begin{equation}\label{pw}
-y_i''+q_iy_i = \mu y_i, \quad i=1,\dots,n,
\end{equation}
or 
\begin{equation}\label{nw}
-y_i''+q_iy_i = -\mu y_i, \quad i= n+1, \dots, K,
\end{equation}
with boundary conditions (\ref{dir}).

Let $\lambda _1^D \leq  \lambda _2^D \leq  \dots$ be the eigenvalues (repeated according to multiplicity) of the system (\ref{dir}) with (\ref{pw}) and (\ref{nw}) for which the eigenvectors are in ${\cal L}^2(G)\cap C^+$. Let $\Lambda ^D_1 < \Lambda ^D_2< \dots $ be the eigenvalues of the system (\ref{dir}) with (\ref{pw}) and (\ref{nw}) not repeated by multiplicity. Denote by $\nu _j^D$ the dimension of the maximal positive (with respect to $[\cdot ,\cdot ]$) subspace of the eigenspace $E_j^D$ to $\Lambda _j^D$.

Observe that if $\mu $ is an eigenvalue of the system (\ref{dir}) with (\ref{pw}) and (\ref{nw}), with multiplicity $\nu $ and eigenspace $E$, then there are precisely $\nu $ indices $i_1, \dots, i_{\nu }$ such that $\mu $ is an eigenvalue of 
\begin{equation}\label{now}
-y_i''+q_iy_i = b_i\mu y_i,
\end{equation}
with boundary conditions (\ref{dir}). In particular, if  
\[ Y^i_j := \left\{\begin{array}{cc}
 0, & j\not= i, \\
y_i, & j=i,
\end{array}\right.\] 
where $j \in \{1, \dots, K\}$, then $Y^{i_1}, \dots, Y^{i_{\nu }}$ are eigenfunctions to (\ref{dir}) with (\ref{pw}) and (\ref{nw}) and  form a basis for $E$, which is orthogonal with respect to both $(\cdot ,\cdot )$ and $[\cdot ,\cdot ]$. Hence, by \cite[Corollary 10.1.4]{Lan}, the maximal $B$-positive subspace of $E$ has dimension $$\nu ^+ = \# \, (\{i_1, \dots, i_{\nu }\}\cap \{1, \dots, n\}).$$ I.e. $\nu ^+$ is the multiplicity of $\mu $ as an eigenvalue of (\ref{dir}) with (\ref{pw}).

Hence $\lambda _j^D$ is the $j$th eigenvalue of (\ref{dir}) with (\ref{pw}), i.e. of (\ref{diff}) with (\ref{dir}) considered only on $G^+$. 

Similarly if we consider the equation
(\ref{diff-1}) with the non-Dirichlet conditions
\begin{equation}\label{nond}
 y_i'(1)= f(1)y_i(1) \quad \text{and} \quad y_i'(0)= f(0)y_i(0), \quad i=1,\dots, K,
\end{equation}
where $f$ is given in (\ref{varform}), then, as in the Dirichlet case, above,
$G$ decomposes into a union of disconnected graphs $G_1, \dots, G_K$. Again, depending on whether the edge has positive or negative weight, we have the equation
\begin{equation}\label{pwnd}
-y_i''+q_iy_i = \mu y_i, \quad i=1,\dots,n,
\end{equation}
or 
\begin{equation}\label{nwnd}
-y_i''+q_iy_i = -\mu y_i, \quad i = n+1, \dots, K,
\end{equation}
with boundary conditions (\ref{nond}).

Let $\lambda _1^N \leq  \lambda _2^N \leq  \dots$ be the eigenvalues (repeated according to multiplicity) of the system (\ref{nond}) with (\ref{pwnd}) and (\ref{nwnd}) for which the eigenvectors are in ${\cal L}^2(G)\cap C^+$. By the same reasoning as above, $\lambda _j^N$ is the $j$th eigenvalue of (\ref{nond}) with (\ref{pwnd}), i.e. of (\ref{diff}) with (\ref{nond}) considered only on $G^+$.  

Thus, from Theorem \ref{mm} and \cite{scbw1} we have that, in ${\cal L}^2(G)\cap C^+$, the eigenvalues of (\ref{diff-1}), (\ref{gjbc1})-(\ref{gjbc2}) are ordered by 
\begin{eqnarray}\label{counting-2}
 \lambda_n^N\le\lambda_n\le\lambda_n^D,\quad n=1,2, \dots.
\end{eqnarray}

The asymptotics for $\lambda _n^N$ and $\lambda _n^D$ are well known, in particular, using the results in \cite{scbw1} for (\ref{diff}) on $G^+$, with (\ref{dir}) and (\ref{nond}) we obtain the following theorem:

\begin{thm}\label{As}
Let $G$ be a compact graph with finitely many nodes. If the boundary value problem (\ref{diff-1}), 
(\ref{gjbc1})-(\ref{gjbc2}) has co-normal (elliptic) boundary conditions, then the eigenvalues in ${\cal L}^2(G)\cap C^+$ obey the asymptotic development
\[\sqrt{\lambda _n} = \frac{n\pi}{\makebox{length}(G^+)} + O(1),\quad \mbox{as}\ n\to\infty.\]
\end{thm} 

By formally replacing $\lambda $ by $-\lambda $ in (\ref{diff}) a similar result is obtained for ${\cal L}^2(G)\cap C^-$.


\end{document}